\documentclass[12pt]{article}

\usepackage{amsmath, epsfig, cite}
\usepackage{amssymb}
\usepackage{amsfonts}
\usepackage{latexsym}
\usepackage{amsthm}

\makeatletter
\renewcommand{\@seccntformat}[1]{{\csname the#1\endcsname}.\hspace{.5em}}
\makeatother

\newtheorem{thm}{Theorem}[section]

\newtheorem{prob}[thm]{Problem}
\newtheorem{cor}[thm]{Corollary}
\newtheorem{conj}[thm]{Conjecture}
\newtheorem{lem}[thm]{Lemma}

\newcommand{\pf}{\noindent{\it Proof.} }

\def\con{{\rm conj}}

\renewcommand{\qed}{\hfill$\Box$\medskip}

\numberwithin{equation}{section}

\begin{document}

\begin{center}
{\Large\bf Some congruences involving central $q$-binomial coefficients}
\end{center}

\vskip 2mm \centerline{Victor J. W. Guo$^1$  and Jiang Zeng$^{2}$}
\begin{center}
{\footnotesize $^1$Department of Mathematics, East China Normal University,\\ Shanghai 200062,
 People's Republic of China\\
{\tt jwguo@math.ecnu.edu.cn,\quad http://math.ecnu.edu.cn/\textasciitilde{jwguo}}\\[10pt]
$^2$Universit\'e de Lyon; Universit\'e Lyon 1; Institut Camille
Jordan, UMR 5208 du CNRS;\\ 43, boulevard du 11 novembre 1918,
F-69622 Villeurbanne Cedex, France\\
{\tt zeng@math.univ-lyon1.fr,\quad
http://math.univ-lyon1.fr/\textasciitilde{zeng}} }
\end{center}

\vskip 0.7cm \noindent{\small{\bf Abstract.}} Motivated by  recent works of Sun and Tauraso,
we prove some variations on the Green-Krammer identity
involving central $q$-binomial coefficients, such as
\begin{align*}
\sum_{k=0}^{n-1}(-1)^kq^{-{k+1\choose 2}}{2k\brack k}_q
\equiv \left(\frac{n}{5}\right) q^{-\lfloor n^4/5\rfloor} \pmod{\Phi_n(q)},
\end{align*}
where $\big(\frac{n}{p}\big)$ is the Legendre   symbol and $\Phi_n(q)$ is the $n$th cyclotomic polynomial.
As consequences, we deduce that
\begin{align*}
\sum_{k=0}^{3^a m-1} q^{k}{2k\brack k}_q
&\equiv 0 \pmod{(1-q^{3^a})/(1-q) },  \\
\sum_{k=0}^{5^a m-1}(-1)^kq^{-{k+1\choose 2}}{2k\brack k}_q
&\equiv 0 \pmod{(1-q^{5^a})/(1-q)},
\end{align*}
for $a,m\geq 1$, the first one being a partial $q$-analogue of the
Strauss-Shallit-Zagier congruence modulo powers of $3$. Several related conjectures are proposed.



\section{Introduction}

The $p$-adic order of several sums involving central binomial coefficients have attracted much  attention.
For example, among other things, Pan and Sun \cite{PS} and Sun and Tauraso \cite{ST,ST2} proved the
following congruences modulo a prime $p$:
\begin{align}
\sum_{k=0}^{p^a-1}{2k\choose k+d}      &\equiv \left(\frac{p^a-|d|}{3}\right) \pmod{p}, \label{eq:st-1}\\
\sum_{k=0}^{p^a-1}(-1)^k{2k\choose k}  &\equiv \left(\frac{p^a}{5}\right) \pmod{p},  \label{eq:st-2}
\end{align}
and, for $p\geq 3$,
\begin{align}
\sum_{k=0}^{p^a-1}{2k\choose k+d}2^{-k}
&\equiv
\begin{cases}
0\pmod p &\text{if $p^a\equiv |d|\pmod 2$,}\\
1\pmod p &\text{if $p^a\equiv |d|+1\pmod 4$,}\\
-1\pmod p &\text{if $p^a\equiv |d|-1\pmod 4$,}
\end{cases}   \label{eq:st-3}
\end{align}
where $\big(\frac{n}{p}\big)$ is the Legendre symbol.
It is well-known that  binomial identities or congruences usually have nice  $q$-analogues (see \cite{Andrews74}).
Recently Tauraso \cite{Tauraso} has noticed that an identity of Greene-Krammer \cite{Greene} can be served as an inspiration
for searching $q$-analogues of some identities in \cite{ST,ST2}, and,
in particular, he has proved the following  generalization  of \eqref{eq:st-1}:
\begin{align}
\sum_{k=0}^{n-1} q^k {2k\brack k+d}_q  \equiv \left(\frac{n-|d|}{3}\right)q^{\frac{3}{2}r(r+1)+|d|(2r+1)}
\pmod{\Phi_n(q)}  \label{eq:k+d}
\end{align}
with $r=\lfloor 2(n-|d|)/3\rfloor$. Here and in what follows
$\Phi_n(q)$ denotes the $n$th cyclotomic polynomial, and
 ${n\brack k}_q$ is
the $q$-binomial coefficient defined by
\[
{n\brack k}_q=
\begin{cases}
\displaystyle\frac{(q;q)_{n}}{(q;q)_{k}(q;q)_{n-k}}, &\text{if $0\leq k\leq n$},\\[5pt]
0,&\text{otherwise,}
\end{cases}
\]
where $(z;q)_n=(1-z)(1-zq)\cdots(1-zq^{n-1})$ is  the  $q$-shifted factorial for $n\geq 0$.

The purpose of this paper is to study  some $q$-versions
 of \eqref{eq:st-1}--\eqref{eq:st-3} as well as some variations of the same flavor as in \cite{Tauraso}.
 For example, from a
$q$-analogue of \eqref{eq:st-2},
we will deduce the following two congruences:
\begin{align}
\sum_{k=0}^{3^a m-1} q^{k}{2k\brack k}_q
&\equiv 0 \pmod{(1-q^{3^a})/(1-q) },  \label{eq:final-1}\\
\sum_{k=0}^{5^a m-1}(-1)^kq^{-{k+1\choose 2}}{2k\brack k}_q
&\equiv 0 \pmod{(1-q^{5^a})/(1-q)}.  \label{eq:final-2}
\end{align}
Note that \eqref{eq:final-1} may be deemed to be a partial $q$-analogue of the
Strauss-Shallit-Zagier congruence \cite{SSZ}:
\begin{align*}
\sum_{k=0}^{3^{a}m-1}{2k\choose k}\equiv 0 \pmod{3^{2a}}.
\end{align*}

The rest of the paper is organized as follows. In Section \ref{sec:q-mod5} we will give a
$q$-analogue of \eqref{eq:st-2} by using  a finite Rogers-Ramanujan identity due to Schur. In Section \ref{sec:cong35}
we will prove \eqref{eq:final-1} and \eqref{eq:final-2}.
Some different $q$-analogues of \eqref{eq:st-3} will be given in Section \ref{sec:2kmod} and
some open problems will be proposed in the last section.

\section{A $q$-analogue of \eqref{eq:st-2}\label{sec:q-mod5}}
It was conjectured by Krammer and proved by Greene \cite{Greene} that
\begin{align}
1+2\sum_{k=1}^{n-1}(-1)^k q^{-{k\choose 2}}{2k-1\brack k}_q
=\begin{cases}
\displaystyle \left(\frac{m}{5}\right)\sqrt{5}, &\text{$n\equiv 0\pmod 5$},\\[10pt]
\displaystyle \left(\frac{n}{5}\right), &\text{otherwise,}
\end{cases} \label{eq:greene}
\end{align}
where $q=e^{2\pi m i/n}$ with $\gcd(m,n)=1$ (see also \cite{Andrews92,CS} for some related results).
If $n=p^a$, then the left-hand side of \eqref{eq:greene} is a $q$-analogue of that of \eqref{eq:st-2}.
However, we cannot deduce the Sun-Tauraso congruence \eqref{eq:st-2} from \eqref{eq:greene}
in the case $n\equiv 0\pmod 5$.
In this section we shall give a new $q$-series identity which is similar to \eqref{eq:greene}
and will imply the Sun-Tauraso congruence \eqref{eq:st-2} completely.


\begin{thm}\label{thm:pimod}
For $n\geq 0$, there holds
\begin{align}
\sum_{k=0}^{n-1}(-1)^kq^{-{k+1\choose 2}}{2k\brack k}_q
\equiv \left(\frac{n}{5}\right) q^{-\lfloor n^4/5\rfloor} \pmod{\Phi_n(q)}. \label{eq:mod5}
\end{align}
In other words,
letting $\omega=e^{2\pi m i/n}$ with $\gcd(m,n)=1$, we have
\begin{align}
&\sum_{k=0}^{n-1}(-1)^k\omega^{-{k+1\choose 2}}{2k\brack k}_\omega
=
\begin{cases}
0,                         &\text{if $n\equiv 0\pmod 5$,} \\
\omega^{-\lfloor n/5\rfloor},   &\text{if $n\equiv 1\pmod 5$,} \\
-\omega^{-\lfloor 3n/5\rfloor}, &\text{if $n\equiv 2\pmod 5$,} \\
-\omega^{-\lfloor 2n/5\rfloor}, &\text{if $n\equiv 3\pmod 5$,} \\
\omega^{-\lfloor 4n/5\rfloor},   &\text{if $n\equiv 4\pmod 5$.}
\end{cases}
\label{eq:pimod5}
\end{align}
\end{thm}

\pf
Since $\omega^k\neq 1$ for $1\leq k\leq n-1$ and $\omega^n=1$, we can write
\begin{align}
\omega^{-{k+1\choose 2}}{2k\brack k}_{\omega}
&=\omega^{-{k+1\choose 2}}\prod_{j=1}^k\frac{1-\omega^{2k+1-j}}{1-\omega^j}\nonumber \\
&=\omega^{-{k+1\choose 2}} (-1)^k \omega^{(3k^2+k)/2}\prod_{j=1}^k\frac{1-\omega^{n-(2k+1-j)}}{1-\omega^j}\nonumber\\
&=(-1)^k\omega^{k^2} {n-k-1\brack k}_{\omega}.\label{eq:lem2.1}
\end{align}
Therefore, we derive from Schur's identity (see, for example, \cite[p.~50]{Andrews}) that
\begin{align}
\sum_{k=0}^{n-1}(-1)^k\omega^{-{k+1\choose 2}}{2k\brack k}_{\omega}
&=\sum_{k=0}^{n-1}\omega^{k^2} {n-k-1\brack k}_{\omega} \nonumber \\
&=\sum_{j=-\infty}^{\infty} (-1)^j \omega^{\frac{j(5j+1)}{2}}{n-1\brack \lfloor\frac{n-5j-1}{2}\rfloor}_\omega.
\label{eq:appschur}
\end{align}
Since $\omega^n=1$, we have
$$
{n-1\brack k}_\omega
=\prod_{i=1}^k\frac{1-\omega^{n-i}}{1-\omega^i}
=\prod_{i=1}^k\frac{1-\omega^{-i}}{1-\omega^i}
=(-1)^k \omega^{-{k+1\choose 2}}
$$
for $0\leq k\leq n-1$, and the result is easily deduced. For example,
if $n=5m$, then there are $2m$ non-zero terms in the right-hand side of \eqref{eq:appschur}.
But the terms indexed $j=-m+2k$ and $j=-m+2k+1$ cancel each other for $k=0,\ldots,m-1$.
\qed

Replacing $q$ by $q^{-1}$, one sees that \eqref{eq:mod5} is equivalent to
\begin{align*}
\sum_{k=0}^{n-1}(-1)^kq^{-{k\choose 2}}{2k\brack k}_q
\equiv \left(\frac{n}{5}\right) q^{\lfloor n^4/5\rfloor} \pmod{\Phi_n(q)}. 
\end{align*}
If $n=p^a$ is  a prime power,  letting $q=1$ in \eqref{eq:mod5},
one immediately gets the Sun-Tauraso congruence \eqref{eq:st-2} by
the formula
\begin{align}\label{lem:q=1}
\Phi_n(1)=
\begin{cases}
p,&\text{if $n=p^a$ is a prime power},\\
1, &\text{otherwise.}
\end{cases}
\end{align}
(Eq.~\eqref{lem:q=1} follows from the identity $q^n-1=\prod_{d|n}\Phi_{d}(q)$ by induction.)

\medskip
\noindent{\it Remark.} The first part of the proof of Theorem \ref{thm:pimod} can be generalized as follows.
We define the $q$-Fibonacci polynomials (see \cite{Carlitz}) by $F_0^q(t)=0$, $F_1^q(t)=1$, and
\begin{align*}
F_n^q(t)=F_{n-1}^q(t)+q^{n-2}tF_{n-2}^q(t),\quad n\geq 2.
\end{align*}
The following is an explicit formula for the $q$-Fibonacci polynomials:
\begin{align}
F_n^q(t)=\sum_{k\geq 0}q^{k^2}{n-k-1\brack k}_q t^k. \label{eq:fibo}
\end{align}
Let $n> d\geq 0$ and $\omega$ be as in Theorem \ref{thm:pimod}. Similarly to \eqref{eq:lem2.1}, we have
\begin{align*}
\omega^{-{k-d\choose 2}}{2k\brack k+d}_{\omega}
&=\omega^{-{k-d\choose 2}} (-1)^{k-d} \omega^{(3k+d+1)(k-d)/2}
\prod_{j=1}^{k-d}\frac{1-\omega^{2n-(2k+1-j)}}{1-\omega^j}\nonumber\\
&=(-1)^{k-d}\omega^{(k+d+1)(k-d)} {2n-k-d-1\brack k-d}_{\omega},
\end{align*}
which yields the following congruence
\begin{align*}
\sum_{k=0}^{n-1}q^{-{k-d\choose 2}}{2k\brack k+d}_{q} t^k
\equiv t^d F_{2(n-d)}^q(-tq^{2d+1}) \pmod{\Phi_n(q)}
\end{align*}
by applying \eqref{eq:fibo}.

\section{Congruences modulo $\Phi_{3^j}(q)$ and $\Phi_{5^j}(q)$ \label{sec:cong35}}
In this section we give a proof of \eqref{eq:final-1} and \eqref{eq:final-2}. It is well-known that
\begin{align*}
\frac{1-q^{p^a}}{1-q}=\prod_{j=1}^a \Phi_{p^j}(q)
\end{align*}
for any prime $p$. We need the following two lemmas.
\begin{lem}\label{lem:n-kk}
For $n\geq 0$, there holds
\begin{align*}
\sum_{k=0}^{n}(-1)^{k} q^{k\choose 2}{n-k\brack k}_q
=(-1)^{n} \left(\frac{n+1}{3}\right)q^{\frac{n(n-1)}{6}}.
\end{align*}
\end{lem}

\begin{lem}\label{lem:root}
Let $m,k,d$ be positive integers, and write $m=ad+b$ and $k=rd+s$,
where $0\leq b,s\leq d-1$. Let $\omega$ be a primitive $d$-th root
of unity.  Then
$$
{m\brack k}_{\omega}={a\choose r}{b\brack s}_\omega.
$$
\end{lem}
\noindent{\it Remark.}
Lemma \ref{lem:n-kk} has appeared in the literature from different origins  (see \cite{Cigler}).
A proof  using mathematical induction
is given in \cite{Tauraso}  and  a multiple extension is proposed  in
\cite{GZ}. Lemma \ref{lem:root} is equivalent to the  $q$-Lucas theorem (see \cite{Olive} and  \cite[Proposition 2.2]{De}).
\medskip

We first establish the following theorem.
\begin{thm}\label{thm:mn-to-35}
Let $m,n\geq 1$. Then
\begin{align}
\sum_{k=0}^{mn-1}q^{k}{2k\brack k}_q
&\equiv \sum_{j=0}^{m-1}{2j\choose j}\sum_{k=0}^{n-1} q^{k}{2k\brack k}_q
\pmod {\Phi_n(q)},  \label{eq:mn-to3}\\
\sum_{k=0}^{mn-1}(-1)^k q^{-{k+1\choose 2}}{2k\brack k}_q
&\equiv \sum_{j=0}^{m-1}(-1)^j{2j\choose j}\sum_{k=0}^{n-1}(-1)^k q^{-{k+1\choose 2}}{2k\brack k}_q
\pmod {\Phi_n(q)}.  \label{eq:mn-to5}
\end{align}
\end{thm}
\pf Let $q=\omega$ be a primitive $n$th root of unity. Then $\omega^n=1$ and
\begin{align}
\sum_{k=jn}^{jn+n-1} \omega^{k}{2k\brack k}_\omega
&=\sum_{k=0}^{n-1} \omega^{k}{2jn+2k\brack jn+k}_\omega, \label{eq:jnn-1}\\
\sum_{k=jn}^{jn+n-1}(-1)^k \omega^{-{k+1\choose 2}}{2k\brack k}_\omega
&=(-1)^{jn}\sum_{k=0}^{n-1}(-1)^k \omega^{-{jn+k+1\choose 2}}{2jn+2k\brack jn+k}_\omega.
\label{eq:jnn-2}
\end{align}
By Lemma \ref{lem:root} we have
$$
{2jn+2k\brack jn+k}_\omega
={2j\choose j}{2k\brack k}_\omega,
$$
which is equal to $0$ if $2k\geq n$. Noticing that
$$
\omega^{-{jn+k+1\choose 2}}=\omega^{-{jn+1\choose 2}}\cdot \omega^{-{k+1\choose 2}}
$$
and
$$
(-1)^{jn}\omega^{-{jn+1\choose 2}}=(-1)^j,
$$
we can write Eqs.~\eqref{eq:jnn-1} and \eqref{eq:jnn-2} as
\begin{align}
\sum_{k=jn}^{jn+n-1}\omega^{k}{2k\brack k}_\omega
&={2j\choose j} \sum_{k=0}^{n-1}\omega^{k}{2k\brack k}_\omega, \label{eq:jnn-11} \\
\sum_{k=jn}^{jn+n-1}(-1)^k \omega^{-{k+1\choose 2}}{2k\brack k}_\omega
&=(-1)^j {2j\choose j}
\sum_{k=0}^{n-1}(-1)^k \omega^{-{k+1\choose 2}}{2k\brack k}_\omega.
\label{eq:jnn-22}
\end{align}
Summing \eqref{eq:jnn-11} and \eqref{eq:jnn-22} over $j$ from $0$ to $m-1$, we complete the proof.
\qed

We now state our main theorem in this section.
\begin{thm} Let $a\geq 1$ and $m\geq 1$. Then
\begin{align}
\sum_{k=0}^{3^a m-1} q^{k}{2k\brack k}_q
&\equiv 0
\pmod{\prod_{j=1}^a \Phi_{3^j}(q)},  \label{eq:phi3} \\
\sum_{k=0}^{5^a m-1}(-1)^kq^{-{k+1\choose 2}}{2k\brack k}_q
&\equiv 0
\pmod{\prod_{j=1}^a \Phi_{5^j}(q)}.   \label{eq:phi5}
\end{align}
\end{thm}
\pf Let $\omega$ be a primitive $n$th
root of unity. Then
$$
\omega^{k}{2k\brack k}_\omega=\omega^{k^2+k}{2k\brack k}_{\omega^{-1}}
=\con\left(\omega^{-k^2-k}{2k\brack k}_{\omega} \right),
$$
where $\con(z)$ denotes the complex conjugate of $z\in\mathbb{C}$.
 From \eqref{eq:lem2.1} we deduce that
$$
\omega^{-k^2-k}{2k\brack k}_{\omega}
=(-1)^k \omega^{{k\choose 2}} {n-k-1\brack k}_{\omega}.
$$
Therefore, by Lemma \ref{lem:n-kk}, we have
\begin{align*}
\sum_{k=0}^{n-1}\omega^{k}{2k\brack k}_\omega
=\con\left(\sum_{k=0}^{n-1}(-1)^k \omega^{{k\choose 2}} {n-k-1\brack k}_{\omega}\right)
=(-1)^{n-1} \left(\frac{n}{3}\right)\omega^{-\frac{(n-1)(n-2)}{6}}.
\end{align*}
This implies that
\begin{align}
\sum_{k=0}^{n-1} q^{k}{2k\brack k}_q \equiv 0 \pmod{\Phi_{n}(q)}
\quad\text{if}\quad 3|n, \label{eq:2kbrack}
\end{align}
which also follows directly from Tauraso's congruence \eqref{eq:k+d}.

Now, letting $n=3^j$ with $1\leq j\leq a$ in \eqref{eq:2kbrack} and
letting  $n=5^j$ with $1\leq j\leq a$ in \eqref{eq:mod5}, we get
\begin{align}
\sum_{k=0}^{3^j-1} q^{k}{2k\brack k}_q &\equiv 0 \pmod{\Phi_{3^j}(q)},  \label{eq:3j} \\
\sum_{k=0}^{5^j-1}(-1)^kq^{-{k+1\choose 2}}{2k\brack k}_q &\equiv 0 \pmod{\Phi_{5^j}(q)}.
\end{align}
Letting $m\to 3^{a-j}m$, $n\to 3^j$ in \eqref{eq:mn-to3} and $m\to 5^{a-j}m$, $n\to 5^j$ in \eqref{eq:mn-to5}
respectively, we obtain
\begin{align*}
\sum_{k=0}^{3^a m-1} q^{k}{2k\brack k}_q
&\equiv 0 \pmod{ \Phi_{3^j}(q)}\quad(1\leq j\leq a),  \\
\sum_{k=0}^{5^a m-1}(-1)^kq^{-{k+1\choose 2}}{2k\brack k}_q
&\equiv 0 \pmod{ \Phi_{5^j}(q)} \quad(1\leq j\leq a).
\end{align*}
Since the cyclotomic polynomials are pairwise relatively prime, we complete the proof.
\qed

We have the following conjecture.

\begin{conj}\label{conj:mod3}
Let $a\geq 1$ and $m\geq 1$. Then
\begin{align*}
\sum_{k=0}^{3^a m-1} q^{k}{2k\brack k}_q  &\equiv 0      \pmod{\prod_{j=1}^a \Phi_{3^j}^2(q)},\\
\sum_{k=0}^{5^{a}-1}(-1)^{k}{2k\choose k} &\equiv 5^{a}  \pmod{5^{a+1}}.
\end{align*}
\end{conj}

We now give a dual form of Theorem \ref{thm:pimod}. The reader is encouraged to compare
it with \cite[Theorem 5.1]{Tauraso}.
\begin{thm}
Let $q=e^{2\pi m i/n}$ with $\gcd(m,n)=1$. Then
\begin{align*}
\sum_{k=0}^{n-1}q^{2k+1}{2k\brack k}_q
=\begin{cases}
\displaystyle \left(\frac{m}{3}\right)i\sqrt{3}, &\text{if $3| n$},\\[10pt]
\displaystyle \left(\frac{n}{3}\right), &\text{otherwise.}
\end{cases}
\end{align*}
\end{thm}

\pf First note that
$$
q^{2k+1}{2k\brack k}_q=q^{(k+1)^2}{2k\brack k}_{q^{-1}}
=\con\left(q^{-(k+1)^2}{2k\brack k}_{q} \right)
$$
and $\Phi_n(q)=0$. From \eqref{eq:lem2.1} we deduce that
$$
q^{-(k+1)^2}{2k\brack k}_{q}
=(-1)^k q^{{k-1\choose 2}-2} {n-k-1\brack k}_{q}.
$$
Therefore,
$$
\sum_{k=0}^{n-1}q^{2k+1}{2k\brack k}_q
=\con\left(\sum_{k=0}^{n-1}(-1)^k q^{{k-1\choose 2}-2} {n-k-1\brack k}_{q}\right).
$$
Since
$$
q^{{k-1\choose 2}-2} {n-k-1\brack k}_{q}
=q^{-n}\left(q^{k+1\choose 2} {n-k\brack k+1}_{q}-q^{k+1\choose 2} {n-k-1\brack k+1}_{q}\right),
$$
by Lemma \ref{lem:n-kk} we have
\begin{align*}
\sum_{k=0}^{n-1}(-1)^k q^{k-1\choose 2} {n-k-1\brack k}_{q}
=(-1)^{n} \left(\left(\frac{n+2}{3}\right)q^{\frac{n(n-5)}{6}}
+\left(\frac{n+1}{3}\right)q^{\frac{n(n-7)}{6}}\right).
\end{align*}
The result then follows easily.
\qed

\begin{cor}
For any positive integer $n$ with $\gcd(n,3)=1$, there holds
\begin{align*}
\sum_{k=0}^{n-1}q^{2k+1}{2k\brack k}_q\equiv \left(\frac{n}{3}\right)\pmod{\Phi_{n}(q)}.
\end{align*}
\end{cor}

For the following  remarkable  congruence of  Sun and Tauraso \cite[(1.1) with $d=0$]{ST}:
\begin{align} \label{eq:newconj2}
\sum_{k=0}^{p^a-1}{2k\choose k}
\equiv \left(\frac{p^a}{3}\right)\pmod{p^2},
\end{align}
we have two interesting $q$-versions to offer:
\begin{conj}
Let $p$ be a prime and $a\geq 1$. Then
\begin{align*}
\sum_{k=0}^{p^a-1}q^{k}{2k\brack k}_q
\equiv \left(\frac{p^a}{3}\right)
q^{\left\lfloor \frac{p^a}{2}-\left(\frac{p^a}{3}\right)\frac{p^a}{6}\right\rfloor
+\left\lfloor\frac{p^a}{3}\right\rfloor p^a}\pmod{\Phi_{p^a}^2(q)},
\end{align*}
and, for $p\neq 3$,
\begin{align*}
\sum_{k=0}^{p^a-1}q^{2k+1}{2k\brack k}_q
\equiv \left(\frac{p^a}{3}\right)
q^{\left(\left\lfloor\frac{p^a+1}{3}\right\rfloor+\left(\frac{p^a}{3}\right)\right)p^a}
\pmod{\Phi_{p^a}^2(q)}.
\end{align*}
\end{conj}

\section{Some $q$-analogues of \eqref{eq:st-3} \label{sec:2kmod}}

To give $q$-analogues of \eqref{eq:st-3}, we need to establish the following $q$-series
identities:

\begin{thm}\label{thm1} Let $n\geq 1$ and $d=0,1,\ldots,n$. Then
\begin{align}
&\sum_{k=0}^{n}(-1)^{n-k}q^{{n-k\choose 2}}{n\brack k}_q{2k\brack k+d}_{q}(-q^{k+1}; q)_{n-k}
=\begin{cases}
\displaystyle q^{\frac{n^2-d^2}{2}}{n\brack \frac{n-d}{2}}_{q^2}, &\text{if $n-d$ is even} \\
\displaystyle 0, &\text{if $n-d$ is odd,}
\end{cases}   \label{eq:qser-1}\\
&\sum_{k=0}^{n}(-1)^{n-k}q^{n-k\choose 2} {n\brack k}_q{2k\brack k+d}_{q}(-q^k; q)_{n-k}
=\begin{cases}
\displaystyle q^{\frac{n^2-d^2}{2}}{n\brack \frac{n-d}{2}}_{q^2}, &\text{if $n-d$ is even} \\
\displaystyle q^{\frac{n^2-d^2-1}{2}}(q^n-1){n-1\brack \frac{n-d-1}{2}}_{q^2}, &\text{if $n-d$ is odd.}
\end{cases}   \label{eq:qser-2}
\end{align}
\end{thm}
\begin{proof}
The $d=0$ case of \eqref{eq:qser-1} was found by Andrews \cite[Theorem 5.5]{Andrews74}.
Both \eqref{eq:qser-1} and \eqref{eq:qser-2} can be proved similarly by using Andrews's
$q$-analogue of Gauss's second theorem
 \cite{Andrews73,Andrews74}:
\begin{align}\label{eq:andrews}
\sum_{k=0}^\infty \frac{(a;q)_k(b;q)_k q^{k+1\choose 2}}{(q;q)_k(abq;q^2)_k}
=\frac{(-q;q)_\infty (aq;q^2)_\infty (bq;q^2)_\infty}{(abq;q^2)_\infty},
\end{align}
where $(z;q)_\infty=\lim_{n\to\infty}(z;q)_n$. We first sketch the proof of \eqref{eq:qser-1}.

Recall that  $(q;q)_{2n}=(q;q^2)_n(q^2;q^2)_n$, $(a;q)_n(-a;q)_n=(a^2;q^2)_n$ and
$$
(a;q)_{n-k}=\frac{(a;q)_n}{(q^{1-n}/a;q)_k}\left(-\frac{q}{a}\right)^k q^{{k\choose 2}-nk}.
$$
Replacing $k$ by $n-k$, we can write the left-hand side of \eqref{eq:qser-1} as
\begin{align*}
&\sum_{k=0}^{n}(-1)^{k} q^{{k\choose 2}}{n\brack k}_q{2n-2k\brack n-k+d}_q(-q^{n-k+1}; q)_{k}\\
&={2n\brack n+d}_q
\sum_{k=0}^n\frac{(q^{-n-d};q)_k(q^{-n+d};q)_kq^{k(k+1)/2}}
{(q;q)_k(q^{-2n+1};q^2)_k}\\
&={2n\brack n+d}_q
\frac{(-q;q)_\infty (q^{-n-d+1};q^2)_\infty (q^{-n+d+1};q^2)_\infty} {(q^{-2n+1};q^2)_\infty}
\quad\text{(by \eqref{eq:andrews})}\\
&=\begin{cases}
\displaystyle q^{\frac{n^2-d^2}{2}}{n\brack \frac{n-d}{2}}_{q^2}, &\text{if $n-d$ is even} \\
\displaystyle 0, &\text{if $n-d$ is odd.}
\end{cases}
\end{align*}
This proves  \eqref{eq:qser-1}.

Observing that
\begin{align*}
\frac{(a;q)_k(b;q)_k q^{k\choose 2}}{(q;q)_k(abq;q^2)_k}
-\frac{(a;q)_k(b;q)_k q^{k+1\choose 2}}{(q;q)_k(abq;q^2)_k}
=\frac{(1-a)(1-b)(aq;q)_{k-1}(bq;q)_{k-1} q^{k\choose 2}}{(1-abq)(q;q)_{k-1}(abq^3;q^2)_{k-1}},
\end{align*}
we derive the following $q$-series identity from \eqref{eq:andrews}:
\begin{align}\label{eq:andrews-2}
\sum_{k=0}^\infty \frac{(a;q)_k(b;q)_k q^{k\choose 2}}{(q;q)_k(abq;q^2)_k}
=\frac{(-q;q)_\infty (aq;q^2)_\infty (bq;q^2)_\infty}{(abq;q^2)_\infty}
+\frac{(-q;q)_\infty (a;q^2)_\infty (b;q^2)_\infty}{(abq;q^2)_\infty}.
\end{align}

Replacing $k$ by $n-k$, we can write the left-hand side of \eqref{eq:qser-2} as
\begin{align*}
&\sum_{k=0}^{n}(-1)^{k} q^{{k\choose 2}}{n\brack k}_q{2n-2k\brack n-k+d}(-q^{n-k}; q)_{k}\\
&={2n\brack n+d}_q
\sum_{k=0}^n\frac{(q^{-n-d};q)_k(q^{-n+d};q)_kq^{k+1\choose 2}}{(q;q)_k(q^{-2n+1};q^2)_k}\frac{(1+q^{n-k})}{(1+q^n)}\\
&={2n\brack n+d}_q
\frac{(-q;q)_\infty (q^{-n-d+1};q^2)_\infty (q^{-n+d+1};q^2)_\infty} {(q^{-2n+1};q^2)_\infty} \\
&\quad{}+\frac{q^n}{1+q^n}{2n\brack n+d}_q
\frac{(-q;q)_\infty (q^{-n-d};q^2)_\infty (q^{-n+d};q^2)_\infty} {(q^{-2n+1};q^2)_\infty}
\quad\text{(by \eqref{eq:andrews} and \eqref{eq:andrews-2})}\\
&=\begin{cases}
\displaystyle q^{\frac{n^2-d^2}{2}}{n\brack \frac{n-d}{2}}_{q^2}, &\text{if $n-d$ is even} \\
\displaystyle q^{\frac{n^2-d^2-1}{2}}(q^n-1){n-1\brack \frac{n-d-1}{2}}_{q^2}, &\text{if $n-d$ is odd.}
\end{cases}
\end{align*}
This proves  \eqref{eq:qser-2}.
\end{proof}

Since $q^n\equiv 1\pmod{\Phi_n(q)}$ and
\begin{align*}
{n-1\brack k}_q=\prod_{j=1}^k\frac{1-q^{n-j}}{1-q^j} &\equiv (-1)^k q^{-{k+1\choose 2}} \pmod{\Phi_n(q)}, \\
{n-1\brack k}_{q^2}=\prod_{j=1}^k\frac{1-q^{2n-2j}}{1-q^{2j}} &\equiv (-1)^k q^{-k(k+1)} \pmod{\Phi_n(q)}, \\
{n-2\brack k}_{q^2}=\prod_{j=1}^k\frac{1-q^{2n-2j-2}}{1-q^{2j}} &\equiv (-1)^k q^{-k(k+3)}
\frac{1-q^{2k+2}}{1-q^2} \pmod{\Phi_n(q)},
\end{align*}
we obtain the following result by substituting $n$ with $n-1$ in \eqref{eq:qser-1} and \eqref{eq:qser-2}.
\begin{cor}Let $n\geq 1$ and $d=0,1,\ldots,n-1$. Then
\begin{align}
&\hskip -3mm
\sum_{k=0}^{n-1}q^{k}{2k\brack k+d}_q(-q^{k+1}; q)_{n-k-1} \nonumber \\
&\equiv
\begin{cases}
\displaystyle 0, &\text{if $n-d$ is even}\\
\displaystyle (-1)^{\frac{n+d-1}{2}}q^{\frac{d(2n-3d)-(n+1)^2}{4}}, &\text{if $n-d$ is odd}
\end{cases}
\pmod{\Phi_n(q)},  \label{eq:2k-first} \\
&\hskip -3mm
\sum_{k=0}^{n-1}q^{k}{2k\brack k+d}_q(-q^{k}; q)_{n-k-1} \nonumber \\
&\equiv
\begin{cases}
\displaystyle (-1)^{\frac{n+d}{2}}q^{\frac{d(2n-3d)-n^2+2d}{4}}\frac{1-q^{n-d}}{1+q}, &\text{if $n-d$ is even}\\
\displaystyle (-1)^{\frac{n+d-1}{2}}q^{\frac{d(2n-3d)-(n+1)^2}{4}}, &\text{if $n-d$ is odd}
\end{cases}
\pmod{\Phi_n(q)}.  \label{eq:2k-second}
\end{align}
\end{cor}

Replacing $q$ by $q^{-1}$ in \eqref{eq:2k-first}, we get
\begin{align*}
&\hskip -3mm
\sum_{k=0}^{n-1}q^{-{k+1\choose 2}}{2k\brack k+d}_q(-q^{k+1}; q)_{n-k-1} \nonumber \\
&\equiv
\begin{cases}
\displaystyle 0, &\text{if $n-d$ is even}\\
\displaystyle (-1)^{\frac{n+d-1}{2}}q^{\frac{1-(n-d)^2}{4}}, &\text{if $n-d$ is odd}
\end{cases}
\pmod{\Phi_n(q)}.
\end{align*}
We also have the following variant of Theorem~\ref{thm1}.
\begin{thm}\label{thm2} Let $n\geq 1$ and $d=0,1,\ldots,n$. Then
\begin{align}
&\sum_{k=0}^{n}(-q)^{n-k}{n\brack k}_q{2k\brack k+d}_{q}(-q^{k+1}; q)_{n-k}
=\begin{cases}
\displaystyle {n\brack \frac{n-d}{2}}_{q^2}, &\text{if $n-d$ is even} \\
\displaystyle (1-q^{2n}){n-1\brack \frac{n-d-1}{2}}_{q^2}, &\text{if $n-d$ is odd.}
\end{cases}  \label{eq:qser-3}\\
&\sum_{k=0}^{n}(-q)^{n-k}{n\brack k}_q{2k\brack k+d}_{q}(-q^k; q)_{n-k}
=\begin{cases}
\displaystyle {n\brack \frac{n-d}{2}}_{q^2}, &\text{if $n-d$ is even} \\
\displaystyle (1-q^n){n-1\brack \frac{n-d-1}{2}}_{q^2}, &\text{if $n-d$ is odd.}
\end{cases}  \label{eq:qser-4}
\end{align}
\end{thm}
\begin{proof} We would only prove \eqref{eq:qser-3}, since the proof of \eqref{eq:qser-4}
is similar. Replacing $q$ by $q^{-1}$ and multiplying by $q^{n^2-d^2+1}$,
one sees that \eqref{eq:qser-3} is equivalent to
the following identity:
\begin{align}
&\hskip -3mm
\sum_{k=0}^{n}(-1)^{n-k}q^{n-k-1\choose 2}{n\brack k}_q{2k\brack k+d}_{q}(-q^{k+1}; q)_{n-k} \nonumber\\
&=\begin{cases}
\displaystyle q^{\frac{n^2-d^2+2}{2}} {n\brack \frac{n-d}{2}}_{q^2}, &\text{if $n-d$ is even} \\
\displaystyle q^{\frac{(n-1)^2-d^2}{2}}(q^{2n}-1){n-1\brack \frac{n-d-1}{2}}_{q^2}, &\text{if $n-d$ is odd.}
\end{cases}  \label{eq:qser-4.0}
\end{align}
Replacing $k$ by $n-k$, we can write the left-hand side of \eqref{eq:qser-4.0} as
\begin{align*}
&\sum_{k=0}^{n}(-1)^{k} q^{{k-1\choose 2}}{n\brack k}_q{2n-2k\brack n-k+d}_q(-q^{n-k+1}; q)_{k}\\
&={2n\brack n+d}_q
\sum_{k=0}^n\frac{(q^{-n-d};q)_k(q^{-n+d};q)_kq^{{k\choose 2}+1}}
{(q;q)_k(q^{-2n+1};q^2)_k}\\
&=q{2n\brack n+d}_q
\frac{(q^{-n-d+1};q^2)_\infty (q^{-n+d+1};q^2)_\infty+(q^{-n-d};q^2)_\infty (q^{-n+d};q^2)_\infty}
 {(q;q^2)_\infty(q^{-2n+1};q^2)_\infty}
\quad\text{(by \eqref{eq:andrews-2})},
\end{align*}
which is equal to the right-hand side of \eqref{eq:qser-4.0}.
\end{proof}

\medskip
\noindent{\it Remark.} Whenever they are discovered, both  Theorem~\ref{thm1} and Theorem~\ref{thm2}
 can be proved by the $q$-Zeilberger algorithm
(see, for example,\cite[p.~113]{Koepf}).
\medskip

As before, we have the following consequences.

\begin{cor}Let $n\geq 1$ and $d=0,1,\ldots,n-1$. Then
\begin{align}
&\hskip -3mm
\sum_{k=0}^{n-1}q^{-k(k+3)/2}{2k\brack k+d}_q(-q^{k+1}; q)_{n-k-1} \nonumber \\
&\equiv
\begin{cases}
\displaystyle (-1)^{\frac{n+d-2}{2}}q^{\frac{5-(n-d+1)^2}{4}}(1-q^{n-d}), &\text{if $n-d$ is even}\\
\displaystyle (-1)^{\frac{n+d-1}{2}}q^{\frac{5-(n-d)^2}{4}}, &\text{if $n-d$ is odd}
\end{cases}
\pmod{\Phi_n(q)},  \label{eq:2k-third}  \\
&\hskip -3mm
\sum_{k=0}^{n-1}q^{-k(k+3)/2}{2k\brack k+d}_q(-q^{k}; q)_{n-k-1} \nonumber \\
&\equiv
\begin{cases}
\displaystyle (-1)^{\frac{n+d-2}{2}}q^{\frac{9-(n-d+1)^2}{4}}\frac{1-q^{n-d}}{1+q}, &\text{if $n-d$ is even}\\
\displaystyle (-1)^{\frac{n+d-1}{2}}q^{\frac{5-(n-d)^2}{4}}, &\text{if $n-d$ is odd}
\end{cases}
\pmod{\Phi_n(q)}.   \label{eq:2k-fourth}
\end{align}
\end{cor}

If we change $q$ to $q^{-1}$, then the congruence \eqref{eq:2k-third} may be rewritten as
\begin{align*}
&\hskip -3mm
\sum_{k=0}^{n-1}q^{2k}{2k\brack k+d}_q(-q^{k+1}; q)_{n-k-1} \nonumber \\
&\equiv
\begin{cases}
\displaystyle (-1)^{\frac{n+d}{2}}q^{\frac{d(2n-3d)-n^2+2d-4}{4}}(1-q^{n-d}), &\text{if $n-d$ is even}\\
\displaystyle (-1)^{\frac{n+d-1}{2}}q^{\frac{d(2n-3d)-(n+1)^2-4}{4}}, &\text{if $n-d$ is odd}
\end{cases}
\pmod{\Phi_n(q)},
\end{align*}
while the congruences \eqref{eq:2k-second} and \eqref{eq:2k-fourth} exchange each other.

\section{Open problems}


Inspired by the $q=1$ case of congruences \eqref{eq:final-1}--\eqref{eq:final-2} and the work of Sun {\cite{Sun},
we would like to make the following conjectures:

\begin{conj}\label{conj-1}
Let $p$ be a prime factor of $4m-1$ with $m\in \mathbb{Z}$ and let $a,n\geq 1$. Then
\begin{align*}
\sum_{k=0}^{p^a n-1}{2k\choose k}m^k\equiv 0 \pmod {p^a}.
\end{align*}
\end{conj}

\begin{conj}\label{conj-2}
Let $m$ be a positive integer. Then
\begin{align*}
\sum_{k=0}^{4m-2}{2k\choose k}m^k\equiv 0 \pmod {(4m-1)},  \\
\sum_{k=0}^{4m}{2k\choose k}(-m)^k\equiv 0 \pmod {(4m+1)}.
\end{align*}
\end{conj}

It is easy to see that Conjecture \ref{conj-1} implies Conjecture \ref{conj-2} but not vice versa.

\begin{conj}Let $a$ be a positive integer. Then
\begin{align*}
\sum_{k=0}^{3^a-1}(-2)^{k}{2k\choose k} &\equiv 3^a \pmod{3^{a+1}}, \\
\sum_{k=0}^{3^a-1}(-5)^{k}{2k\choose k} &\equiv 2\cdot3^a \pmod{3^{a+1}}, \\
\sum_{k=0}^{7^a-1}(-5)^{k}{2k\choose k} &\equiv 7^a \pmod{7^{a+1}}.
\end{align*}
\end{conj}

\begin{conj}\label{conj-4}
Let $m$ be a positive integer. If $4m-1$ is a prime and $m\neq 1$, then
\begin{align*}
\sum_{k=0}^{(4m-1)^a-1}{2k\choose k}m^k\equiv (4m-1)^a \pmod {(4m-1)^{a+1}}.
\end{align*}
If $4m+1$ is a prime, then
\begin{align*}
\sum_{k=0}^{(4m+1)^a-1}{2k\choose k}(-m)^k\equiv (4m+1)^a \pmod {(4m+1)^{a+1}}.
\end{align*}
\end{conj}

Conversely, we make the following conjecture, which gives a sufficient condition for
whether $4m-1$ or $4m+1$ is a prime.
We have checked the cases $m\leq 1500$ via Maple, not finding any counter examples.

\begin{conj}
Let $m$ be a positive integer. If $m\neq 30$ and
\begin{align*}
\sum_{k=0}^{4m-2}{2k\choose k}m^k\equiv 4m-1 \pmod {(4m-1)^2},
\end{align*}
then $4m-1$ is a prime. If
\begin{align*}
\sum_{k=0}^{4m}{2k\choose k}(-m)^k\equiv 4m+1 \pmod {(4m+1)^2},
\end{align*}
then $4m+1$ is a prime.
\end{conj}

The following conjecture looks a little different but seems also very challenging.
\begin{conj}\label{conj:fin}
Let $a$ and $n$ be positive integers. Then
\begin{align}
\sum_{k=0}^{5^{a}n-1}{4k\choose 2k}{2k\choose k}^2
&\equiv 0 \pmod{5^{a}},  \label{eq:fin-conj}\\
\sum_{k=0}^{5^{a}-1}{4k\choose 2k}{2k\choose k}^2
&\equiv (-1)^a5^a \pmod{5^{a+1}}.  \nonumber
\end{align}
\end{conj}

\noindent{\it Remark.} Recently, Pan and Sun \cite{PS2} have
confirmed the first congruence in Conjecture \ref{conj:mod3} and Sun \cite{Sun2} has proved
Conjectures \ref{conj-1}--\ref{conj-4}
(naturally including the second congruence in Conjecture \ref{conj:mod3}).

\begin{prob}Are there any $q$-analogues of Conjectures \ref{conj-1}--\ref{conj:fin}{\rm?}
\end{prob}

\vskip 2mm \noindent{\bf Acknowledgments.} The authors thank the anonymous referee
for helpful comments on this paper. The first author was sponsored by
Shanghai Educational Development Foundation under the Chenguang
Project (\#2007CG29), Shanghai Rising-Star Program (\#09QA1401700),
Shanghai Leading Academic Discipline Project (\#B407),
and the National Science Foundation of China (\#10801054).
The second author was supported by the project MIRA 2008 of R\'egion Rh\^one-Alpes.

\end{document}